\documentclass[letterpaper, 10 pt, conference]{ieeeconf}
\IEEEoverridecommandlockouts                             % This command is only
\usepackage[english]{babel}

\usepackage{verbatim}
\usepackage{amsmath}
\usepackage{amssymb}
\usepackage{graphicx}

\usepackage[labelformat=simple]{subcaption}

\DeclareCaptionLabelSeparator{periodspace}{.\quad}
\usepackage{hyperref}

\usepackage[usenames,dvipsnames]{pstricks}
 \usepackage{epsfig}
 \usepackage{pst-grad} % For gradients
 \usepackage{pst-plot} % For axes

\title{\LARGE \bf Pattern Formation with a Compartmental Lateral Inhibition System}

\author{\authorblockN{Ana S. Rufino Ferreira, Justin Hsia, Murat Arcak, Michel Maharbiz and Adam Arkin}
\authorblockA{
ana, jhsia, arcak, maharbiz@eecs.berkeley.edu and aparkin@lbl.gov}
}

\newtheorem{theorem}{Theorem}[section]

\newtheorem{lemma}[theorem]{Lemma}

\newtheorem{assumption}[theorem]{Assumption}
\newtheorem{claim}[theorem]{Claim}

\newcommand{\R}{\mathbb{R}}

\renewcommand{\baselinestretch}{1}

\begin{document}
\maketitle
\thispagestyle{empty}
\pagestyle{empty}

\begin{abstract}
We propose a compartmental lateral inhibition system that generates contrasting patterns of gene expression between neighboring compartments.  The system consists of a set of compartments interconnected by channels.  Each compartment contains a colony of cells that produce diffusible molecules to be detected by the neighboring colony, and each cell is equipped with an inhibitory circuit that reduces its production when the detected signal is stronger.  We develop a technique to analyze the steady-state patterns emerging from this lateral inhibition system and apply it to a specific implementation. The analysis shows that the proposed system indeed exhibits contrasting patterns  within realistic parameter ranges.
\vspace{-1mm}
\end{abstract}

\section{Introduction}

Multicellular developmental processes rely on spatial patterning to initiate differentiation \cite{gilbert10,wolpert11}.
Commonly-studied methods of pattern formation include diffusion-driven instability \cite{turing1952,meinhardt1982,hsia2012}, gradient or density detection \cite{basu2005,liu2011}, locally-synchronized oscillators \cite{danino10}, and lateral inhibition \cite{kunisch94,collier96,sprinzak10}.
Lateral inhibition is a mechanism where cell-to-cell signaling induces neighboring cells to compete and diverge into sharply contrasting fates, enabling developmental processes such as segmentation or boundary formation \cite{meinhardt2000}.
The best-known example of lateral inhibition is the Notch pathway in Metazoans where membrane bound Delta ligands bind to the Notch receptors on the neighboring cells. This binding releases the Notch intracellular domain in the neighbors, which then inhibits their Delta ligand production \cite{muskavitch1994,collier96,sprinzak11,arcak13}.
Recent discoveries have shown that lateral inhibition is not limited to complex organisms: a contact-dependent inhibition (CDI) system has been identified in \textit{E. coli} where delivery via membrane-bound proteins of the C-terminus of the gene \textit{cdiA} causes down regulation of metabolism \cite{aoki2009,aoki2010,webb2013}. Despite the vigorous research on elucidating natural pathways such as Notch and CDI, a \textit{synthetic} lateral inhibition system for pattern formation has not been developed.\\

\begin{figure}[ht]
\centering
\includegraphics[page=2,width=0.1888\textwidth]{pstricksfigures3index}%0.18
\includegraphics[page=1,width=0.288\textwidth]{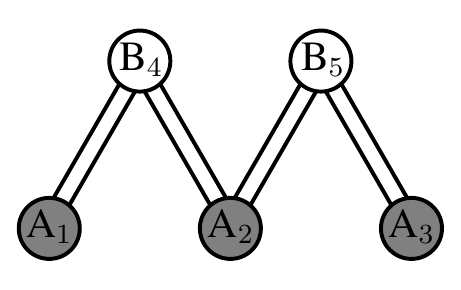}%0.285
\vspace{-1mm}\caption{Compartmental lateral inhibition system with cells of type $A$ and $B$, where contrasting patterns between neighboring compartments emerge. In each compartment $A_i$ ($B_i$) we place a colony with cells of type $A$ ($B$) that communicate through channels. Each cell type can only detect signaling molecules produced by the other type, preventing auto-inhibition.}
\label{fig:spatialcontactnetwork}
\vspace{-2mm}
\end{figure}

In this paper, we propose a \textit{compartmental} lateral inhibition system that is able to spontaneously generate contrasting patterns between neighboring compartments. Our system consists of a set of compartments interconnected by channels as in Figure \ref{fig:spatialcontactnetwork}. In each compartment, we place a colony of cells that produce diffusible molecules to be detected by the neighboring colony. We equip each cell with an inhibitory circuit that reacts to the detected signal, \textit{i.e.}, the more diffusible molecules are detected in one compartment, the less production in that colony.
To prevent auto-inhibition, we use two orthogonal diffusible quorum sensing molecules \cite{collins06} and design two inhibitory circuits each of which detects only one type of molecule and produces the other type.
In the examples of Figure \ref{fig:spatialcontactnetwork}, cells of type $A$ produce a diffusible molecule $X$ that is only detectable by cells of type $B$, and cells of type $B$ produce a diffusible molecule $Y$ which is only detectable by $A$.
\vspace{-3mm}

To derive conditions under which this system will exhibit contrasting patterns, we define the cell network as a graph where each compartment corresponds to a vertex. The diffusion of molecules between two compartments occurs through the channels and is represented by the graph edges. We model the diffusion with a compartmental model, and represent the compartment-to-compartment communication by the Laplacian matrix of the weighted graph. The edge weights depend on the distance between the compartments and the diffusivity of the quorum sensing molecules. We then use the graph-theoretic notion of \textit{equitable partition} to ascertain the existence of contrasting steady-state patterns. Equitable partitions reduce the steady-state analysis to finding the fixed-points of a scalar map, and each fixed-point represents a steady-state where all the compartments of the same type have the same final value. We also show that the slope of the scalar map at each fixed-point provides a stability condition for the respective steady-states. Finally, we propose and model a synthetic circuit with cells of type $A$ and $B$, which is currently under implementation, and apply our analysis to show that it is capable of patterning.\\

\vspace{-1mm}Graph theoretical results have been used in the analysis of patterning by lateral inhibition in our recent work \cite{arcak13, ferreira13}. 
However, these references addressed a contact inhibition model for networks of identical cells, whereas the present paper allows two cell types which is critical for avoiding auto-inhibition in practice.\\

\vspace{-1mm}Reaction-diffusion mechanisms have been widely used in the past to achieve spatial pattern generation with synthetic systems; mostly relying on one-way communication achieved through either the external spatio-temporal manipulation of the cell's chemical environment \cite{cohen09,sohka09,lucchetta05}, the precise positioning of cells containing different gene networks which secrete or respond to diffusible signals \cite{basu2005,basu04}, or the interplay between cell growth and gene expression \cite{payne13}. A two-way communication mechanism using orthogonal quorum sensing systems has been employed to demonstrate a predator-prey system in \cite{balagadde08}. Unlike these results, this paper achieves spatial patterning by lateral inhibition by using orthogonal quorum sensing systems and by positioning colonies of cells inside compartments that are connected by  channels.

\section{An Analytical Test for Patterning}\label{sec:mainresult}
\subsection{Composing a Compartmental Lateral Inhibition Model}
We propose a network of $N_A$ compartments of type $A$ and $N_B$ compartments of type $B$ that communicate through diffusible molecules. 
Each cell of type $A$ produces diffusible species $X$, and only cells of type $B$ are equipped with a receiver species that binds to $X$ and forms a receiver complex. Similarly, the diffusible species $Y$ is produced by cells of type $B$ and detected by cells of type $A$. We represent the dynamics in each cell type with three modules: the transmitter module where species $X$ (or $Y$) are produced and released; the receiver module where $Y$ (or $X$) is detected, and an inhibitory module which inhibits the transmitter activity in the presence of the receiver complex.\\

\begin{figure}[ht]
\centering\vspace{-4.9mm}
  	\includegraphics[page=2,width=0.6\linewidth]{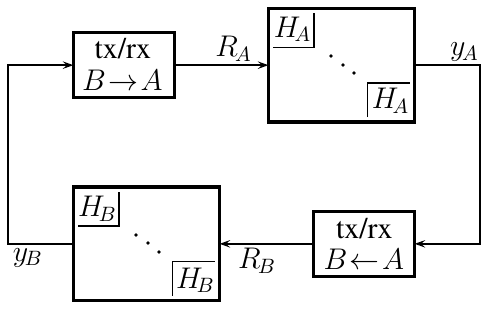}
\vspace{0mm}  \caption{Cell network with two types of compartments $A$ and $B$ communicating through diffusion. For each type of diffusible species, the transceiver includes the dynamics of the transmitter modules on the sender compartments, the receiver modules on the detection compartments, and the diffusion process.}
\label{fig:diffNetDiag}
\end{figure}

\vspace{-0.8mm}To facilitate the analysis, we separate the transmitter module of $A$ and receiver module of $B$, and merge them into a ``transceiver" for the diffusible species $X$, which also includes the diffusion process. Similarly, the transceiver block of $Y$ is composed by the transmitter module of $B$ and the receiver module of $A$. The cell network is represented in Figure \ref{fig:diffNetDiag}.
Each compartment is represented with a block labeled $H_A$ or $H_B$, corresponding to the inhibitory circuit of types $A$ and $B$, respectively. The concentration of the auto-inducer for the production of $X$ (respectively, $Y$) is denoted by $y_A$ ($y_B$), and $R_A$ ($R_B$) is the concentration of the receiver complex, result from the binding of $Y$ ($X$) to the receiver protein.

\vspace{-2mm}The transceiver blocks incorporate diffusion in an ordinary differential equation compartmental model that describes the concentrations of the diffusible species at each compartment.
We define an undirected graph $\mathcal{G}=\mathcal{G}(V,E)$ where each element of the set of vertices $V$ represents one  compartment, and each edge $(i,j)\in E$ represents a channel between compartments $i$ and $j$. For each edge $(i,j)\in E$ we define a weight $d_{ij}=d_{ji}$ (and $d_{ij}=0$ if compartments $i$ and $j$ are not connected). The constant $d_{ij}$ is proportional to the diffusivity of the species and inversely proportional to the square of the distance between compartments $i$ and $j$. We define the weighted Laplacian of the graph to be:\par\nobreak\vspace{-1mm}
{\small\begin{equation}\label{eq:laplacian}
\{L\}_{ij}=\left\{\begin{array}{cc}
-\sum_{j=1}^N d_{ij} & \text{if } i=j\vspace{1mm}\\
d_{ij} & \text{if } i\neq j.
\end{array}\right.
\end{equation}}%

The dynamical model of the transceiver tx/rx for $X$ is then represented by:\par\nobreak\vspace{-4mm}
{\small\begin{equation}\label{eq:transceiverDyn}
\text{tx/rx}_{A\rightarrow B}\hspace{-0.5mm}:\left\{\hspace{-1mm}\begin{array}{ccl}
\left[\hspace{-0.5mm}\begin{array}{c}\dot{X}_A \\ \dot{X}_B\end{array}\hspace{-0.5mm}\right] 
&\hspace{-2mm}{=}\hspace{-2mm}&
\left[\hspace{-0.5mm}\begin{array}{c}\Gamma_X (X_A,y_A) \\ \Phi_X(X_B,R_B)\end{array}\hspace{-0.5mm}\right]
+L\left[\hspace{-0.5mm}\begin{array}{c} X_A\\X_B\end{array}\hspace{-0.5mm}\right]\vspace{1mm}\\
\dot{R}_B 
&\hspace{-2mm}{=}\hspace{-2mm}& \Psi_X(X_B,R_B),
\end{array}
\right.\vspace{-0.5mm}
\end{equation}}%
where $X_A\in\R^{N_A}_{\geq 0}$ represents the concentration of species $X$ in compartments $A$ due to production, $X_B\in\R^{N_B}_{\geq 0}$ the concentration of species $X$ at compartment $B$  due to diffusion, and $R_B\in\R^{N_B}_{\geq 0}$ the concentration of complexes at compartment $B$ formed by the binding of species $X$ with a receiver molecule. The functions $\Gamma_X(\cdot,\cdot)\in\R^{N_A}_{\geq 0}$, $\Phi_X(\cdot,\cdot)\in\R^{N_B}_{\geq 0}$, and $\Psi_X(\cdot,\cdot)\in\R^{N_B}_{\geq 0}$ are concatenations of the decoupled elements $\gamma^i_X(X_A^i,u^i)\,{\in}\,\R_{\geq 0}$, $i\,{=}\,1,...,N_A$, $\phi^j_X(X_B^j,R_B^j)\,{\in}\,\R_{\geq 0}$ and $\psi^j_X(X_B^j,R_B^j)\,{\in}\,\R_{\geq 0}$, $j\hspace{0.4mm}{=}\hspace{0.3mm}1,...,N_B$, and assumed to be continuously differentiable. The transceiver $\text{tx/rx}_{B\rightarrow A}$ for $Y$ is defined similarly, by changing $X$ to $Y$ and switching indices $A$ with $B$ in \eqref{eq:transceiverDyn}.\\

\begin{assumption}\label{ass:differentiabilityTR}
For each constant input $y_A^*\in\R^{N_A}_{\geq 0}$ (and $y_B^*\in\R^{N_B}_{\geq 0}$), the system \eqref{eq:transceiverDyn} has a globally asymptotically stable steady-state $(X_A^*,X_B^*,R_B^*)$, which is an hyperbolic equilibrium, \textit{i.e.}, the Jacobian has no eigenvalues on the imaginary axis. Furthermore, there exist positive and increasing functions $T^{\text{tx/rx}}_{AB}:\R^{N_A}_{\geq 0}\rightarrow\R^{N_B}_{\geq 0}$ and $T^{\text{tx/rx}}_{BA}:\R^{N_B}_{\geq 0}\rightarrow\R^{N_A}_{\geq 0}$ such that
\vspace{-2mm}\begin{equation}\label{eq:TfuncTR}
R_B^*:=T^{\text{tx/rx}}_{AB}(y^*_A),\ \ \ \text{ and }\ \ \   R_A^*:=T^{\text{tx/rx}}_{BA}(y^*_B).\vspace{-2.3mm}
\end{equation}\hfill$\blacksquare$%\Box
\end{assumption}
The increasing property of these maps is meaningful, since a higher input of the auto-inducer leads to more production and, thus, more detection on the receiver side.\\

Next, we represent the blocks $H_k^i$, $i=1,...,N$ of type $k\in\{A,B\}$ with models of the form:
\vspace{-0.3mm}{\begin{equation}\label{eq:cellDyn}
H_{k}^i:\left\{\begin{array}{rcl}
\dot{x}_i & = & f_k(x_i,u_i)\vspace{0.5mm}\\
y_i & = & h_k(x_i),
\end{array}
\right.\vspace{-0.6mm}
\end{equation}}%
where $x_i\in\R^n_{\geq 0}$ describes the vector of reactant concentrations in compartment $i$, $y_i\in\R_{\geq 0}$ is the output of compartment $i$ (in this context, the concentration of an auto-inducer synthase), and $u_i\in\R_{\geq 0}$ is the input of compartment $i$ (the concentration of the receiver complex). We denote $x_k\,{=}\,[x_1^T,...,x_{N_k}^T]^T\in\R^{nN_k}_{\geq 0}$, $u_k\,{=}\,[u_1,...,u_{N_k}]^T\in\R^{N_k}_{\geq 0}$, and $y_k\,{=}\,[y_1,...,y_{N_k}]^T\in\R^{N_k}_{\geq 0}$, $k\in\{A,B\}$.\\

We assume that $f_k(\cdot,\cdot)$ and $h_k(\cdot)$ are continuously differentiable and further satisfy the following properties:
\begin{assumption}\label{ass:differentiability}
For $k\in\{A,B\}$ and each constant input $u^*\in\R_{\geq 0}$, the system
\eqref{eq:cellDyn} has a globally asymptotically stable steady-state
\vspace{-3mm}\begin{equation}
x^*:=S_k(u^*),\vspace{-1mm}
\end{equation}
which is an hyperbolic equilibrium.
Furthermore, the maps $S_k:\R_{\geq 0}\rightarrow\R_{\geq 0}^n$ and
 $T_k:\R_{\geq 0}^n\rightarrow\R_{\geq 0}$, defined as:
\vspace{-0.5mm}\begin{equation}\label{eq:Tfunc}
T_k(\cdot):=h_k(S_k(\cdot)),\vspace{-1mm}
\end{equation}
are continuously differentiable, and $T_k(\cdot)$ is a positive, bounded and decreasing function.\hfill$\blacksquare$%\Box
\end{assumption}
\vspace{2mm}
The decreasing property of $T_k(\cdot)$ is consistent with the lateral inhibition feature, since a higher input in one cell leads to lower output values.

\subsection{When do Contrasting Patterns Emerge?}
We now present a method to find steady-state patterns for the system defined by \eqref{eq:cellDyn}-\eqref{eq:transceiverDyn}. Let $z_A\in\R^{N_A}_{\geq 0}$ and $z_B\in\R^{N_B}_{\geq 0}$ be a steady-state for $y_A$ and $y_B$, respectively. Then, $z_A$ and $z_B$ must satisfy the following:
\begin{equation}\label{eq:bigSSeq}
\left\{\begin{array}{rcl}
z_A &=& \mathbf{T}_A(T^{\text{tx/rx}}_{BA}(\mathbf{T}_B(T^{\text{tx/rx}}_{AB}(z_A))))\vspace{1.5mm}\\
z_B &=& \mathbf{T}_B(T^{\text{tx/rx}}_{AB}(\mathbf{T}_A(T^{\text{tx/rx}}_{BA}(z_B))))
\end{array}\right.\end{equation}
where\par\nobreak\vspace{-5mm}
{\small\begin{eqnarray*}
&\mathbf{T}_{A}(u_A)=[T_A(u_A^1),..., T_A(u_A^{N_A})]^T:\R^{N_A}_{\geq 0}\rightarrow\R^{N_A}_{\geq 0},\\
&\mathbf{T}_{B}(u_B)=[T_B(u_B^1),..., T_B(u_B^{N_B})]^T:\R^{N_B}_{\geq 0}\rightarrow\R^{N_B}_{\geq 0}.
\end{eqnarray*}}%
Given Assumptions \ref{ass:differentiabilityTR} and \ref{ass:differentiability}, a steady-state for $y_A$ and $y_B$ is sufficient to conclude the existence of a steady state for the full system defined by \eqref{eq:transceiverDyn}-\eqref{eq:cellDyn}. Our goal is to determine when $z_A$ and $z_B$ exhibit sharply contrasting values, indicating an on/off pattern.

To reduce the dimension of the maps defined in \eqref{eq:bigSSeq}, we use the notion of \textit{equitable partition} from graph theory \cite[section 9.3]{godsil01}. For a weighted and undirected graph $\mathcal{G}(V,E)$, with a Laplacian matrix $L$ as defined in \eqref{eq:laplacian}, a partition of the vertex set $V$ into classes $O_1,...,O_r$ is said to be \textit{equitable} if there exists $\overline{d}_{ij}$ $i,j=1,...,r$, such that
\begin{equation}\label{eq:equitable}
\sum_{v\in O_j}d_{uv}=\overline{d}_{ij}\ \ \ \forall u\in O_i,\ i\neq j.
\end{equation}
This means that the sum of the edge weights from a vertex in a class $O_i$ into all the vertices in a class $O_j$ ($i\neq j$) is invariant of the choice of the vertex in class $O_i$. We let the \textit{quotient Laplacian} $\overline{L}\in\R^{r\times r}$ be formed by the off-diagonal entries $\overline{d}_{ij}$, and $\left\{\overline{L}\right\}_{ii}=\{L\}_{ii}=-\sum_{j=1,j\neq i}^{r}\overline{d}_{ij}$.\\

\begin{assumption}\label{ass:equitable}
The partition of the compartments $V$ into the classes $O_A$ and $O_B$ of type $A$ and $B$, respectively, is equitable.
\hfill$\blacksquare$
\end{assumption}
\vspace{2mm}
This assumption implies that the total incoming edge weight of the species $X$ is the same for all the compartments of type $B$, and the total incoming edge weight of the species $Y$ is the same for all the compartments of type $A$. For example, the network on the left side of Figure \ref{fig:spatialcontactnetwork} is equitable with respect to the classes $O_A$ and $O_B$ if $d_{13}+d_{14}=d_{23}+d_{24}$ and $d_{13}+d_{23}=d_{14}+d_{24}$, which implies $d_{13}=d_{24}$ and $d_{23}=d_{14}$. Since the edge weights $d_{ij}$ are inversely proportional to the square of the distance, this means that opposite channels must have the same length, thus exhibiting a parallelogram geometry.\\

Assumption \ref{ass:equitable} allows us to search for solutions to \eqref{eq:bigSSeq} where the compartments of the same type have the same steady-state, \textit{i.e.}, 
\begin{equation}\label{eq:smallSS}
z=[\overline{z}_A,...,\overline{z}_A,\overline{z}_B,...,\overline{z}_B]^T= [\overline{z}_A\mathbf{1}_{N_A}^T,\overline{z}_B\mathbf{1}^T_{N_B}]^T
\end{equation}
where $\overline{z}_A\in\R_{\geq 0}$ and $\overline{z}_B\in\R_{\geq 0}$. This means that the transceiver input-output maps become decoupled and $R_B^*=\mathbf{T}_{AB}(\overline{z}_A)$, where $\mathbf{T}_{AB}(\overline{z}_A)=T_{AB}(\overline{z}_A)\mathbf{1}_{N_B}$, with $T_{AB}:\R_{\geq 0}\rightarrow\R_{\geq 0}$; and similarly for $T_{BA}^{\text{tx/rx}}(\cdot)$ with the map $T_{BA}:\R_{\geq 0}\rightarrow\R_{\geq 0}$. Note that the diffusion coefficients are implicit in the maps $T_{AB}(\cdot)$ and $T_{BA}(\cdot)$.\\

The steady-states \eqref{eq:smallSS} must satisfy the following reduced system of equations:
\begin{equation}
\left\{\begin{array}{rcl}\label{eq:smallSSeq}
\overline{z}_A&=&T_A(T_{BA}(T_B(T_{AB}(\overline{z}_A))))\triangleq \overline{T}_A(\overline{z}_A)\vspace{1.5mm}\\
\overline{z}_B&=&T_B(T_{AB}(T_A(T_{BA}(\overline{z}_B))))\triangleq \overline{T}_B(\overline{z}_B)
\end{array}
\right.,\end{equation}
where $\overline{T}_A(\cdot):\R_{\geq 0}\rightarrow\R_{\geq 0}$ and $\overline{T}_B(\cdot):\R_{\geq 0}\rightarrow\R_{\geq 0}$. The solutions of the scalar equations in \eqref{eq:smallSSeq} are solutions of the coupled system of $N_A$ (and $N_B$) equations in \eqref{eq:bigSSeq}. Furthermore, it is sufficient to study the solution of one of the equations in \eqref{eq:smallSSeq}: if $\tilde{z}_A$ is a solution to the top equation, then $\tilde{z}_B\triangleq T_B(T_{AB}(\tilde{z}_A))$ is a solution to the second equation. The derivative of these two functions at the fixed points $\tilde{z}_A$ and $\tilde{z}_B\triangleq T_B(T_{AB}(\tilde{z}_A))$ is the same and given by\par\nobreak
\vspace{-3mm}{\small
\begin{equation}\hspace{-3mm}
\left.\frac{d\overline{T}_A}{dz_A}\right|_{\tilde{z}_A}\hspace{-3.9mm}=\hspace{-0.5mm}T'_{AB}(\tilde{z}_A)T'_B\hspace{-1mm}\left(T_{AB}(\tilde{z}_A)\right)\hspace{-0.5mm}T'_{BA}(\tilde{z}_B)T'_A\hspace{-1mm}\left(T_{BA}(\tilde{z}_B)\right)\hspace{-0.6mm}=\hspace{-1.2mm}\left.\frac{d\overline{T}_B}{dz_B}\right|_{\tilde{z}_B}\hspace{-3mm},\vspace{-0.9mm}
\end{equation}}%
where $T'(z_k)\triangleq\left.\frac{dT}{dz}\right|_{z=z_k}$.

\vspace{-2mm}
\begin{figure}[ht]
\hspace{-1mm}
\begin{subfigure}[b]{.5\linewidth}
\centering%
  	\includegraphics[page=2,width=0.95\linewidth]{grapherPlots}
\vspace{-5mm}\caption{}
\label{fig:TATBfunca}
\end{subfigure}%
\begin{subfigure}[b]{.5\linewidth}
\centering%
   	\includegraphics[page=1,width=0.95\linewidth]{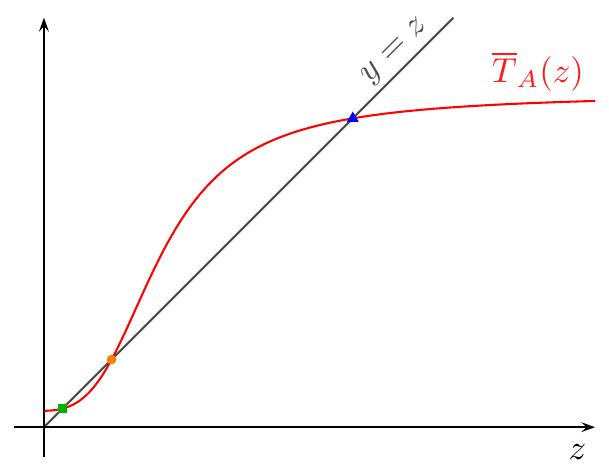}
\vspace{-4mm}\caption{}
\label{fig:TATBfuncb}
\end{subfigure}
\caption{Typical shapes of input-output maps $\overline{T}_A(\cdot)$ and $\overline{T}_B(\cdot)$: (a) In this case, the unique pair of fixed points (orange circles) is near-homogenous and no contrasting patterns emerge; (b) In this case, there exist three pairs of fixed points (orange circle, green square, and blue triangle), and the two extra solutions represent contrasting steady-state patterns.}\label{fig:TATBfuncs1pts}
\end{figure}

From Assumptions \ref{ass:differentiability} and \ref{ass:differentiabilityTR}, the input-output transfer maps $T_A(\cdot)$ and $T_B(\cdot)$ are decreasing and bounded, whereas $T_{BA}(\cdot)$ and $T_{AB}(\cdot)$ are increasing. Since $\overline{T}_A(\cdot)$ and $\overline{T}_B(\cdot)$ in \eqref{eq:smallSSeq} are compositions of these four maps, they are positive, increasing and bounded functions. Figure \ref{fig:TATBfuncs1pts} illustrates typical shapes of the input-output maps $\overline{T}_A(\cdot)$ and $\overline{T}_B(\cdot)$. In \ref{fig:TATBfunca} there exists only one solution pair (orange circles). This is a near-homogeneous steady-state, where the discrepancy between $\tilde{z}_A$ and $\tilde{z}_B$ is due only to nonidentical $\overline{T}_A(\cdot)$ and $\overline{T}_B(\cdot)$. In Figure \ref{fig:TATBfuncb} we see three fixed points: the middle solution pair (near-homogenous steady-state), the largest fixed point in $\overline{T}_A(\cdot)$ (blue triangle), and the smallest fixed point in $\overline{T}_A(\cdot)$ (green square). The latter two have a corresponding opposite fixed-point pair in $\overline{T}_B(\cdot)$, and therefore represent a contrasting steady-state pattern.\\

It is important to note that a contrasting pattern emerges when the near-homogenous steady-state has a slope larger than $1$ as in Figure \ref{fig:TATBfuncb}, that is: 
\begin{equation}\label{eq:existenceCond}
T'_{AB}(\tilde{z}_A)T'_B\left(T_{AB}(\tilde{z}_A)\right)T'_{BA}(\tilde{z}_B)T'_A\left(T_{BA}(\tilde{z}_B)\right)>1.
\end{equation}
Indeed, due to the boundedness and strictly increasing properties of the map $\overline{T}_A(\cdot)$, there must exist at least two other fixed point pairs of \eqref{eq:smallSSeq}, $(z^*_A,z^*_B{\triangleq}T_B(T_{AB}(z^*_A)))$ and $(z^{**}_A,z^{**}_B)$ for which
\vspace{-1mm}
\begin{equation}\label{eq:otherSS}
(z^*_A{>}\tilde{z}_A\ \ \text{and}\ \ z^*_B{<}\tilde{z}_B) \ \ \ (z^{**}_A{<}\tilde{z}_A\ \ \text{and}\ \ z^{**}_B{>}\tilde{z}_B).
\end{equation}
In the next section, we show that \eqref{eq:existenceCond} implies that the near-homogenous steady-state becomes unstable, setting the stage for contrasting patterns to emerge. Thus, \eqref{eq:existenceCond} provides a parameter tuning principle and is instrumental in characterizing the parameter ranges for patterning in Section \ref{sec:example}.

\section{Convergence to Contrasting Patterns}\label{sec:stability}

To analyze convergence to the steady-state patterns in \eqref{eq:smallSSeq}, we employ monotonicity assumptions. A \textit{monotone} system is one that preserves a partial ordering of the initial conditions as the solutions evolve in time, and a partial ordering is defined with respect to a \textit{positivity cone} in the Euclidean space that is closed, convex, pointed ($K\cap (-K)=\{0\}$), and has nonempty interior. In such a cone, $x\preceq\hat{x}$ means $\hat{x}-x\in K$. Given the positivity cones $K^U$, $K^Y$, $K^X$ for the input, output, and state spaces, the system $\dot{x}=f(x,u)$, $y=h(x)$ is said to be \textit{monotone} if $x(0)\preceq \hat{x}(0)$ and $u(t)\preceq\hat{u}(t)$ for all $t\geq 0$ imply that the resulting solutions satisfy $x(t)\preceq \hat{x}(t)$ for all $t\geq 0$, and the output map is such that $x\preceq\hat{x}$ implies $h(x)\preceq h(\hat{x})$ \cite{angeli03}.
\begin{assumption}\label{ass:monotonicityTR}
The system tx/rx$_{A\rightarrow B}$ in \eqref{eq:transceiverDyn} is monotone with respect to $K^U{=}\R_{\geq 0}^{N_A}$, $K^Y{=}\R_{\geq 0}^{N_B}$, and $K^X{=}\R_{\geq 0}^{N+N_B}$. Similarly tx/rx$_{B\rightarrow A}$ is monotone with respect to $K^U{=}\R_{\geq 0}^{N_B}$, $K^Y{=}\R_{\geq 0}^{N_A}$, and $K^X{=}\R_{\geq 0}^{N+N_A}$.\hfill$\blacksquare$
\end{assumption}
\begin{assumption}\label{ass:monotonicity}
The systems $H_A$ and $H_B$ in \eqref{eq:cellDyn} are monotone with respect to $K^U{=}{-}K^Y{=}\R_{\geq 0}$, and $K^X{=}K$, where $K$ is some positivity cone in $\R$.\hfill$\blacksquare$
\end{assumption}
These monotonicity assumptions are consistent with Assumptions \ref{ass:differentiabilityTR} and \ref{ass:differentiability}, as they imply the increasing property of the input-output maps $T^{\text{tx/rx}}_{BA}(\cdot)$ and $T^{\text{tx/rx}}_{AB}(\cdot)$ and the decreasing behavior of $T_A(\cdot)$ and $T_B(\cdot)$.
We now state a stability result for solutions restricted to the steady-state solutions described by \eqref{eq:smallSSeq}.

\vspace{2mm}
\begin{theorem}\label{thm:stabilityFULL}
Consider the network \eqref{eq:transceiverDyn}-\eqref{eq:cellDyn} and suppose Assumptions \ref{ass:differentiabilityTR}, \ref{ass:differentiability}, \ref{ass:monotonicityTR} and \ref{ass:monotonicity} hold. Let the partition of the compartments into the classes $O_A$ and $O_B$ be equitable. Then the steady-state described by  \eqref{eq:smallSSeq} is asymptotically stable if
\vspace{-1mm}\begin{equation}\label{eq:stabilityCondFULL}\vspace{-1mm}
T'_{AB}(\tilde{z}_A)T'_B\left(T_{AB}(\tilde{z}_A)\right)T'_{BA}(\tilde{z}_B)T'_A\left(T_{BA}(\tilde{z}_B)\right)<1,
\end{equation}
and unstable if \eqref{eq:existenceCond} holds.\hfill$\blacksquare$
\end{theorem}
See the Appendix for a proof of this Theorem.

\section{Synthetic Lateral Inhibition Circuit}\label{sec:example}

We propose a lateral inhibition circuit with two types of compartments as described above. The diffusible species are two acyl-homoserine lactones (AHL), namely C10HSL and 3OC6HSL, while the two receiver proteins are LuxR$\mbox{-}$G2E$\mbox{-}$R67M and LuxR, respectively. This choice guarantees that the AHL/LuxR pairs interact orthogonally with each other \cite{collins06}. To keep the notation used in the previous section, we denote C10HSL by $X$, 3OC6HSL by $Y$, and the complexes LuxR$\mbox{-}$G2E$\mbox{-}$R67M$\mbox{-}$C10HSL by $R_B$, and LuxR$\mbox{-}$3OC6HSL by $R_A$.\\

\begin{figure}[ht]
\vspace{0mm}\includegraphics[width=1\linewidth]{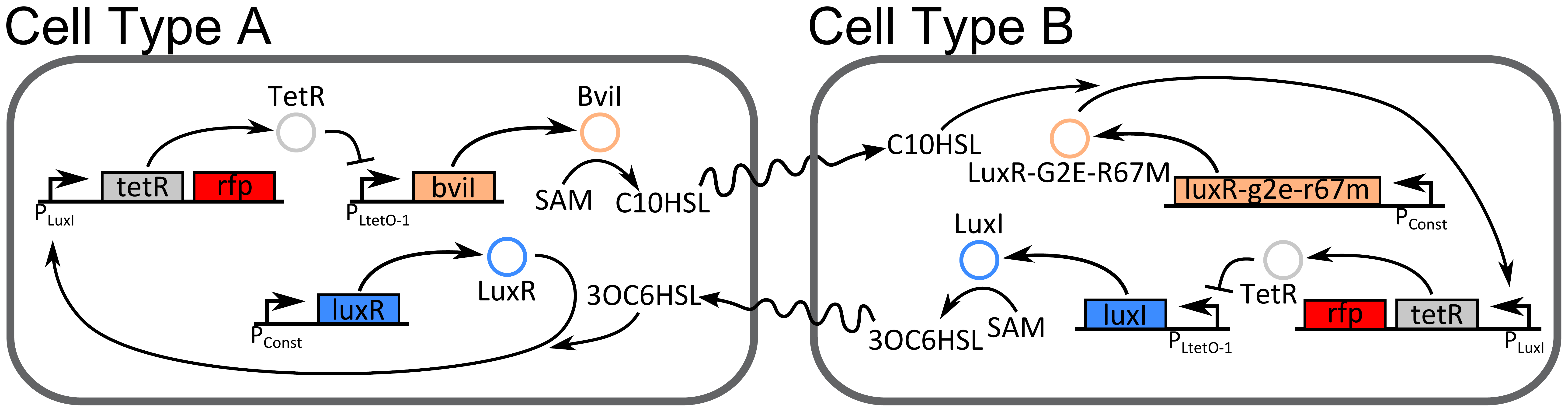}
\caption{Diagram of the synthetic lateral inhibition circuit under implementation, using two orthogonal AHL/LuxR pairs: C10HSL/LuxR$\mbox{-}$G2E$\mbox{-}$R67M and 3OC6HSL/LuxR.}
\label{fig:microfluidicCDI}
\vspace{-3mm}
\end{figure}

In Figure \ref{fig:microfluidicCDI} we represent the synthetic circuit for each cell of type $A$ (left) and $B$ (right). We use the \textit{luxI}/\textit{luxR} (and \textit{bviI}/\textit{luxR$\mbox{-}$g2e$\mbox{-}$r67m}) genes as auto-inducer synthase and receptor, the auto-inducer \textit{luxI} (\textit{bviI}), which is transcribed by P$_{\textit{LtetO}\mbox{-}1}$%\cite{lutz97}
, translates LuxI (BviI) which is responsible for the production of $X$ or $Y$. The receptor proteins, variants of LuxR, in each compartment, detect and bind to the $X$ and $Y$ received, forming the complexes $R_B$ and $R_A$, respectively. The $R_B$ ($R_A$) complex induces the production of the protein TetR and inhibition occurs when TetR represses the promoters P$_{\textit{LtetO}\mbox{-}1}$, thus inhibiting the production of LuxI (BviI). We use red fluorescence protein (RFP) as reporters for each compartment, which are induced by $R_A$ (or $R_B$).

\vspace{-1mm}

We study the following equations for this network, grouped into the transceiver blocks and the inhibitory cell circuits as defined in the Section \ref{sec:mainresult}. The model for the inhibitory circuit of cell type $A$ is:\par\nobreak\vspace{-4mm}
{\small\begin{equation}\label{eq:HA}
H_A^i:\hspace{-1mm}\left\{\hspace{-2.5mm}\begin{array}{rcl}
\frac{d}{dt}m_{T_Y}^i \hspace{-1mm}&\hspace{-2.5mm}=\hspace{-2mm}&\hspace{-1mm} V_{\text{P}_{\text{LuxI}}}N_{\text{P}_{\text{LuxI}}}C \left(\frac{1}{1+
(K_{RA}/R_A^i)^{n_{RA}}}+\ell_{\text{P}_{\text{LuxI}}}\hspace{-0.5mm}\right)+\\
&&\hphantom{12312312312312312312312}-\gamma_{m_T}m_{T_Y}^i\\
\frac{d}{dt}p_{T_Y}^i \hspace{-1mm}&\hspace{-2.5mm}=\hspace{-2mm}&\hspace{-1mm} \epsilon_T m_{T_Y}^i - \gamma_T p_{T_Y}^i\\
\frac{d}{dt}m_{I_X}^i \hspace{-1mm}&\hspace{-2.5mm}=\hspace{-2mm}&\hspace{-1mm} V_{\text{P}_{\text{LtetO-}1}}N_{\text{P}_{\text{LtetO-}1}}C \left(\hspace{-0.5mm}\frac{1}{1+
(p_{T_Y}^i/K_{T})^{n_{T}}}+\ell_{\text{P}_{\text{LtetO-}1}}\hspace{-1mm}\right)+\\
&&\hphantom{12312312312312312312312}-\gamma_{mI_X}m_{I_X}^i\\
\frac{d}{dt}p_{I_X}^i \hspace{-1mm}&\hspace{-2.5mm}=\hspace{-2mm}&\hspace{-1mm} \epsilon_{I_X} m_{I_X}^i - \gamma_{I_X} p_{I_X}^i
\end{array}\right.
\end{equation}}%
where $m_k$ is the $mRNA_k$ concentration and $p_k$ the protein concentration (with subscripts $T$ for TetR and $I$ for LuxI);  $\gamma_k$ the degradation rate; $\epsilon_k$ the translation rate; $V_k$ the transcriptional velocity rate; $N_k$ the copy number; $C$ the concentration (in $M$) of a single molecule in a cell; $\ell_k$ the leakage rate; $K_k$ the dissociation constant; while $n_{RA}$ and $n_T$ represent the cooperativity.\\

For the dynamics of the transceiver of $X$, we consider $X_A^i$, $i=1,...,N_A$ to be the concentration of species $X$ at compartment $i$ of type $A$, and  $X_B^j$, $j=1,...,N_B$ the concentration of species $X$ at compartment $j$ of type $B$. Let $[X^T,R_B^T]$ be the transceiver state, with $X=[X_A^T, X_B^T]^T=[X_A^1,...,X_A^{N_A},X_B^1,...,X_B^{N_B}]^T$ and $R_B=[R_B^1,...,R_B^{N_B}]^T$. The transceiver dynamics are:\par\nobreak\vspace{-3mm}
{\small\vspace{-1mm}\begin{equation}\label{eq:TAB}
\text{tx/rx}_{A\rightarrow B}\hspace{-1mm}:\hspace{-1mm}\left\{\hspace{-2mm}\begin{array}{rcl}
\frac{d}{dt}X_A^i &\hspace{-2mm}=\hspace{-2mm}& \nu p_{I_X}^i - \gamma_{X}X_A^i + L_iX\\
\frac{d}{dt}X_B^j &\hspace{-2mm}=\hspace{-2mm}& -k_{on}X_B^j(p_{R_X} -R_B^j) + k_{off}R_B^j+\\
&&\hphantom{12312312312312}-\gamma_{X}X_B^j+ L_{j+N_A}X\\
\frac{d}{dt}R_B^j &\hspace{-2mm}=\hspace{-2mm}& k_{on}X_B^j(p_{R_X} - R_B^j)-k_{off}R_B^j,\\
\end{array}\right.%\\\vphantom{\sum} \ \ i=1,...,N_A\ \ \ \ j=1,...,N_B\nonumber
\end{equation}}%
for $i=1,...,N_A$, $j=1,...,N_B$, where $L_i$ corresponds to the row $i$ of the Laplacian matrix, $p_{R_k}$ is the constitutive concentration of total LuxR (bound and unbound), $k_{on}/k_{off}$ are the binding rates, and $\nu$ is the generation rate of AHL. The dynamics for the inhibitory circuit of cell type $B$ and for the transceiver ${\text{tx/rx}}_{B\rightarrow A}$ are obtained similarly, by changing the indices appropriately.\\

Next, we analyze the range of parameters where patterning occurs. To analyze the steady-states of the network above, note that both $H_A$ and $\text{tx/rx}_{A\rightarrow B}$ ($H_B$ and $\text{tx/rx}_{B\rightarrow A}$) meet the Assumptions in \ref{ass:differentiability}, \ref{ass:monotonicity} and \ref{ass:differentiabilityTR}, \ref{ass:monotonicityTR}, respectively. From \eqref{eq:HA}, for each constant input $R_A^{i*}$, there is only one steady-state solution $(m_{T_Y}^{i*},p_{T_Y}^{i*},m_{I_X}^{i*},p_{I_X}^{i*})$, which is a globally asymptotically stable hyperbolic equilibrium, due to the lower triangular structure of \eqref{eq:HA} with bounded nonlinearities. Furthermore, the static input-output map is decreasing:
\vspace{-3mm}\begin{equation}\footnotesize
T_{A}^{i}(R_A^{i*}) = K_1\hspace{-1mm}
\left(\hspace{-0.5mm}
\frac{1}{1{+}\hspace{-0.5mm}\left(\hspace{-0.5mm}
\frac{K_2}{K_T}\hspace{-0.5mm}\left(\hspace{-0.5mm}
\frac{1}{1+(K_{RA}/R_A^{i*})^{n_{RA}}}{+}\ell_{\text{P}_{\text{LuxI}}}\right)
\hspace{-0.5mm}\right)^{n_T}} {+} \ell_{\text{P}_{\text{LtetO-}1}}\hspace{-1mm}\right)\hspace{-0.5mm},
\vspace{-1mm}\end{equation}
where\par
{\vspace{-3.5mm}\footnotesize\begin{equation*}
K_1{=}\frac{\epsilon_{I_X}}{\gamma_{I_X}}\frac{V_{\text{P}_{\text{LtetO-}1}}N_{\text{P}_{\text{LtetO-}1}}C}{\gamma_{m_{I_X}}} \ \text{ and } \ K_2{=}\frac{\epsilon_T}{\gamma_T}\frac{V_{\text{P}_{\text{LuxI}}}N_{\text{P}_{\text{LuxI}}}C}{\gamma_{m_T}}\ \  \text{(M)}.
\end{equation*}}%
The subsystem is monotone with respect to $K^U\hspace{-0.5mm}{=}{-}K^Y\hspace{-0.5mm}{=}\hspace{0.2mm}\R_{\geq 0}$, $K{=}\{x{\,\in\,}\R^4|\, x_1{\geq}0,\,x_2{\geq}0,\,x_3{\leq}0,\,x_4{\leq}0\}$ \cite[Lemma 4]{arcak13}.\\

As for the transceiver $\text{tx/rx}_{A\rightarrow B}$ in \eqref{eq:TAB}, we see that in steady-state, for a constant input $p_{I_X}^*\in\R^{N_A}$, the dynamic equations for $R_B$ become zero, which implies that the first terms of the dynamical equations for $X_B$ are also zero. Therefore, due to the linearity of the remainder terms, there exists a unique solution for $[X_A^{*T}, X_B^{*T}]^T$:\par\nobreak\vspace{-0.5mm}
{\small\begin{equation}\label{eq:SStxrx1}
\left[\begin{array}{c}X_A^*\\X_B^*\end{array}\right] = (-L+\gamma_X I_{N})^{-1}\left[\begin{array}{c}\nu p_{I_X}^*\\ 0_{N_B}\end{array}\right].
\end{equation}}%
The inverse of $(-L+\gamma_X I_{N})$ exists since $-L$ is a positive semidefinite matrix (property of Laplacian matrices). The single solution for the steady-state of $R_B^i$ is given by\par\nobreak\vspace{-0.5mm}
{\small\begin{equation}\label{eq:SStxrx2}
R_B^{i*}=\frac{p_{R_X}}{1+\frac{k_{off}}{k_{on}}\frac{1}{X_B^{i*}}},\vspace{-0.5mm}
\end{equation}}%
where $X_B^{i*}$ is as in \eqref{eq:SStxrx1}. Note that the static input-output map $T_{AB}^{tx/rx}(p_{I_X}^{i*})$ is positive and increasing, because $(-L+\gamma_X I_{N})$ is a positive definite matrix with nonpositive off-diagonal elements, and thus its inverse is a positive matrix (\textit{i.e.}, all elements are positive) \cite[Theorem 6.2.3]{berman94}. Finally, to conclude that these steady-states are asymptotically stable and hyperbolic, we write the Jacobian of the transceiver as:\par\nobreak\vspace{-3mm}
{\small\begin{equation}
J=\hspace{-0.5mm}\left[\begin{array}{c|c}
L-\gamma_X I_{N} & \begin{array}{c}{0}\\0\end{array}\\\hline 
\begin{array}{cc}0 & 0\end{array} & 0\end{array}\right]\hspace{-1mm}+\hspace{-1mm}
\left[\begin{array}{c|cc} 0 & 0 & 0\\\hline 0 & -D_{R_B} & D_{X_B}\\ 0 & D_{R_B} & -D_{X_B} \end{array}\right],\vspace{0.5mm}
\end{equation}}%
where $D_{R_B}$ and $D_{X_B}$ are diagonal matrices with elements $\{D_{R_B}\}_{ii}=k_{on}(p_{R_X}{-}R_B^{i*})$ and 
$\{D_{X_B}\}_{ii}=k_{on}X_B^{i*}+k_{off}$, $i=1,...,N_B$. The matrix $J$ has negative diagonal terms and nonnegative off-diagonal terms, and there exists a $D$ such that the column sum of $DJD^{-1}$ are all negative for all states in the nonnegative orthant\footnote{\vspace{-0.5mm}choose $D=diag\{\underbrace{1,...,1}_{N \text{times}},\underbrace{k,...,k}_{N_B \text{times}}\}$, with $1<k<1+\frac{\gamma_X}{k_{on}p_{R_X}}$}. Note that this implies that the \textit{matrix measure} of $DJD^{-1}$ with respect to the one-norm is negative \cite[Chapter 2]{desoer09}, and $\mu_D(J)=\mu_1(DJD^{-1}){<}0$. %, for all $[X^T,R_B^T]^T\in\R_{\geq 0}^{N+N_B}$. 
This is a \textit{contraction} property with respect to the weighted one-norm; therefore, for each constant input, the steady-state is globally asymptotically stable \cite{sontag10}. Moreover, it is an hyperbolic equilibrium since $Re\{\lambda_k(J)\}{\leq}\mu(J){<}0$ \cite{desoer09}. The transceiver is monotone with respect to the cones in Assumption \ref{ass:monotonicityTR} since the Jacobian off-diagonal terms are all positive and the dependence on the input variable $p_{I_X}$ is positive \cite{angeli03}.\\

To find stable steady-state patterns where all the compartments of the same type have the same final value, let the network be an equitable graph $\mathcal{G}$ with respect to the compartment types. The transceiver input-output map decouples into the scalar maps,
\begin{equation}
T_{AB}(\tilde{z}_A) = \frac{1}{1+\frac{k_{off}}{k_{on}}\frac{\gamma_X(\gamma_X+\overline{d_{AB}}+\overline{d_{BA}})}{\overline{d_{BA}}\nu}\frac{1}{\tilde{z}_A}},
\end{equation}
where $\overline{d_{AB}}$ and $\overline{d_{BA}}$ are as in \eqref{eq:equitable}. As discussed in the previous section, we look for the steady-states that are fixed points of $\overline{T}_A(\cdot)$ and $\overline{T}_B(\cdot)$.\\

The reaction parameters used for the analysis are displayed in Table \ref{table:parameters} in the Appendix, and are similar to the parameters suggested in \cite{hsia2012}. We assume that the two orthogonal types of AHL have similar induction and binding reception parameters, and thus consider both cell types to have the same parameter values. In this particular case, the maps $\overline{T}_A(\cdot)$ and $\overline{T}_B(\cdot)$ are identical, and when there exist three fixed-points as in Figure \ref{fig:TATBfuncb}, the middle solution pair is the same for $A$ and $B$ (\textit{i.e.}, $\tilde{z}_A=\tilde{z}_B$). The slope of these maps at the fixed points depends on the edge weights $d_{ij}$ and constitutive concentration of total LuxR $p_{R_i}$, which are tunable parameters. As discussed next, $d_{ij}$ can be tuned by changing the channels' length, and $p_{R_i}$ can be tuned by changing the strength of the constitutive promoter.\\

When each compartment is a square of side $w$, and the channel connecting the compartments be of length $l_{ij}$ and width $w$, the edge weight is, by \cite{dreij11}:
\begin{equation}
d_{ij}=\frac{D_{\text{AHL}}}{l_{ij}w}=k\frac{D_{\text{AHL}}}{l_{ij}^2}.
\vspace{-1mm}\end{equation}
Here we let the width be a factor $k$ of the length, \textit{i.e.} $w=l/k$. In the laboratory, we intend to fill the channel and compartments with agar and pipette one colony in each compartment. As the agar solidifies, a thin layer of water is formed on its surface. The AHL diffusion occurs on the agar surface. Although the cells remain on the agar surface, the AHL diffusion occurs through the agar as well, but we assume this to be negligible in comparison with the diffusion on the surface. We consider the diffusivity coefficient for AHL in water at $25^{\circ}C$ \cite{stewart_b03}: $D_{\text{AHL},25^{^\circ}\hspace{-0.5mm}C}=4.9\times 10^{-10} m^2/s$.\\

\begin{figure}[ht]
\centering
\includegraphics[width=0.79\linewidth]{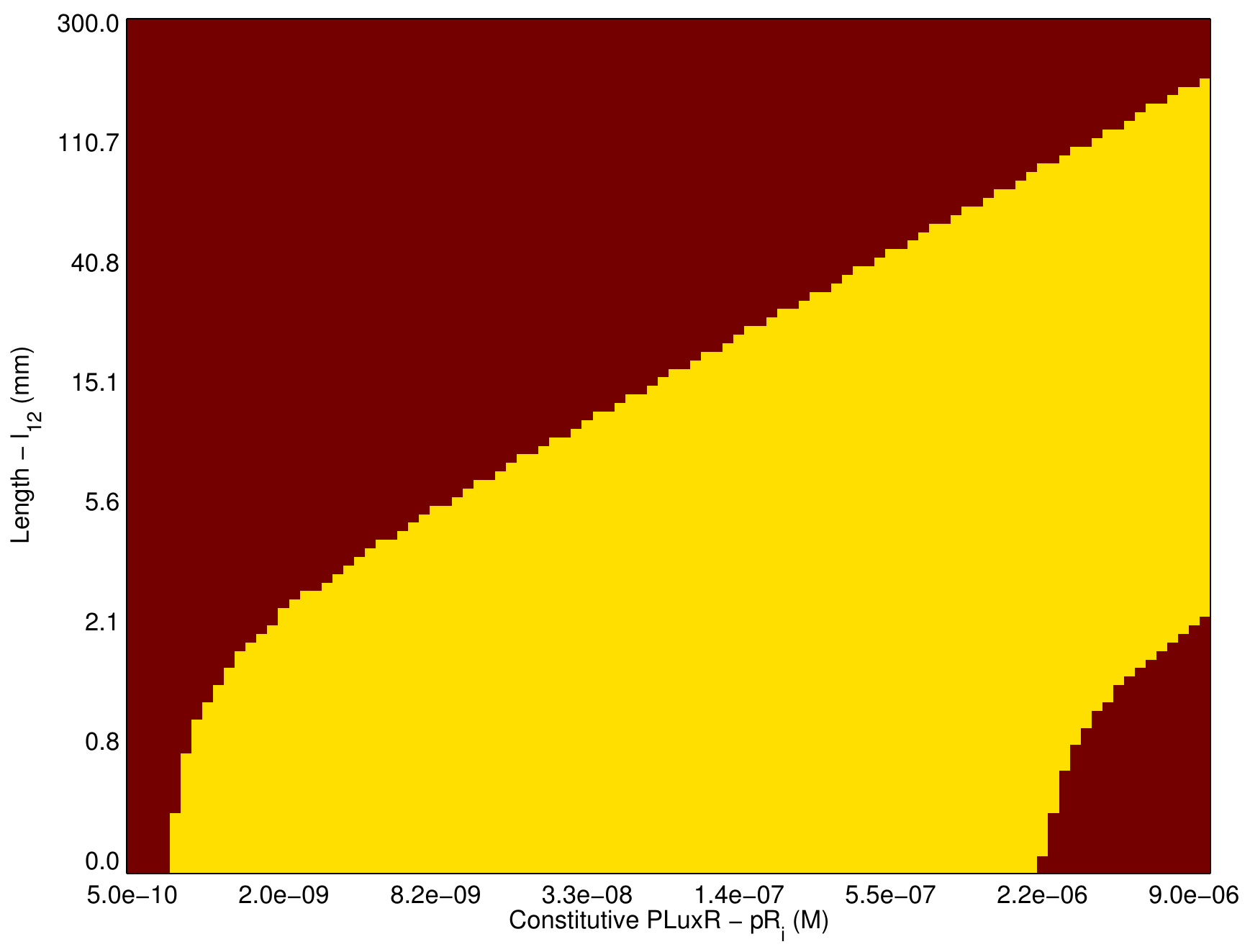}
\caption{Patterning (yellow) vs. Non-Patterning (dark red) region, for varying $p_{R_i}$ and $l_{12}$.}
\label{fig:patterning_PRIvsL12}\vspace{-4mm}
\end{figure}

As an illustration of the patterning condition \eqref{eq:existenceCond}, consider now two compartments connected by one channel, one compartment of type $A$ and the other of type $B$. We assess the slope of the scalar input-output maps by varying the channel length $l_{12}$ and the constitutive concentration of $p_{R_i}$. Figure \ref{fig:patterning_PRIvsL12} maps the regions over the pairs $(p_{R_i},l_{12})$ where contrasting patterns emerge. We obtain patterning within a wide range of realistic values of $p_{R_i}$. At the extreme values, if the concentration of $p_{R_i}$ is too low, the detection ability of each cell is affected, which leads to a low concentration of the complex AHL$\mbox{-}$LuxR, and since no cell is being inhibited (fluorescence reporters are low) contrasting patterning does not occur. When $p_{R_i}$ is too high, the cells are too sensitive to the reception of any leakage AHL, and therefore are inhibited (fluorescence reporters high) and no contrasting patterning occurs.

\begin{figure}[ht]
\centering
\includegraphics[width=0.66\linewidth]{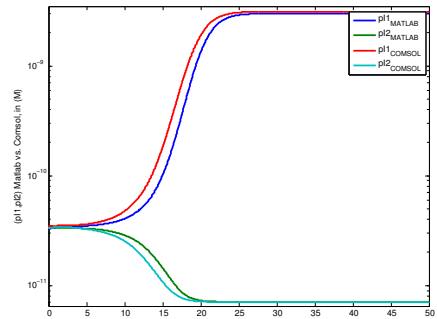}
\caption{Result comparison between ODE model, in MATLAB, and PDE model, in COMSOL (with $p_{R_i}=5\times10^{-7}$ and $l_{12}=500\mu m$). Note that both models converge to the same steady-state, with a similar time constant (${\sim}22$hrs vs. ${\sim}19$hrs, respectively).}
\label{fig:matvscom}\vspace{-4mm}
\end{figure}

There is also a limit on the length of the channel for the emergence of contrasting patterns (Figure \ref{fig:patterning_PRIvsL12}). In implementation, we expect a stricter limit on the length of the channel since the compartmental model does not account for degradation of AHL along the channels. For validation, we have implemented the compartment network in COMSOL, a finite element analysis, solver and simulation software for multi-physics applications, which allows for coupled systems of partial differential equations (PDEs) with complex geometry. In COMSOL, we define the geometry of the channel and the compartments, and only allow AHL to diffuse through the channel. For the values of $p_{R_i}$ studied, we have seen a cap on patterning for lengths no larger than $5$ mm. Due to degradation along the channel, only a small portion of the AHL actually reaches the opposite compartment. Although the ODE model does not account for this, for shorter channels ($\leq 3 mm$), we compute a degradation correction factor for the ODE compartmental model that compensates for the extra degradation along the channel. In these regimes, we obtain an accurate steady-state and dynamical match between the ODE model and the PDE COMSOL model, see Figure \ref{fig:matvscom}.

\section{Conclusions}
In addition to the current effort to implement this design in the laboratory, several analytical problems remain for future research: we will explore the stochastic patterning behavior of these networks, due to the inherent stochastic nature of biochemical reactions and diffusion; and will also investigate a perturbation analysis for the emerging patterns in networks where small variations to the channels' lengths exist, resulting in quasi-equitable graphs.

\section*{Acknowledgment}
This research was supported in part by the NIH National Institute of General Medical Sciences grant 1R01GM109460-01, and by the National Science Foundation grant NSF ECCS-1101876.

\section*{Appendix}
To prove Theorem \ref{thm:stabilityFULL}, we first use the result that the compartmental network  \eqref{eq:transceiverDyn}-\eqref{eq:cellDyn} is monotone. 
\begin{lemma}\label{lemma:monotonicityNetwork}
If monotonicity Assumptions \ref{ass:monotonicityTR} and \ref{ass:monotonicity} hold, then the network \eqref{eq:transceiverDyn}-\eqref{eq:cellDyn} is monotone.\hfill$\blacksquare$
\end{lemma}
We skip this derivation due to space constraints. The main idea of the proof follows similarly to  \cite[Theorem 3]{arcak13}, we can represent the network as a unitary positive feedback interconnection of a monotone system where the inputs and outputs are ordered with respect to the same positivity cone. Note the network is a cascade of an ``anti-monotone" system ($H_A$ composed with tx/rx$_{A\rightarrow B}$) with another ``anti-monotone'' system ($H_B$ composed with tx/rx$_{B\rightarrow A}$), thus the composite system is monotone with the same input and output ordering, $K^U=K^Y=\R_{\geq 0}^{N_A}$ and $K^X=K^{N_A}{\times}\R_{\leq 0}^{N+N_B}{\times}\{-K\}^{N_B}{\times}\R_{\geq 0}^{N+N_A}$.\\

Since the network is monotone, we know from \cite[Lemma 6.4]{angeli04b} that the linearized system around the steady-state is also monotone with respect to the same positivity cones. Furthermore, \cite[Theorem 2]{enciso05} shows that for a linear system $\dot{x}=Ax+Bu$ and $y=Cx$ that is monotone with respect to the cones $K^U{=}K^Y$, $K^X$, and Hurwitz matrix A, the following equivalence holds: $A+BC$ is Hurwitz if and only if $-(I+CA^{-1}B)$ is Hurwitz. Therefore, we can prove stability of the positive feedback monotone system from the ``dc-gain'' of the open loop system.\\

\textit{Proof of Theorem \ref{thm:stabilityFULL}:}
The linearization of the full network \eqref{eq:transceiverDyn}-\eqref{eq:cellDyn} about the steady state is given by:
\begin{equation}\label{eq:bigjacobianFULL}
\arraycolsep=1.4pt\def\arraystretch{1.4}
%\left[\begin{array}{c}\dot{X}_A\\\dot{X}_{AB}\\\dot{X}_B\\\dot{X}_{BA}\end{array}\right]=
\hspace{-1mm}\left[\footnotesize\hspace{-1mm}\begin{array}{cccc}
A_A{\otimes}I_{N_{\hspace{-0.5mm}A}} & 0 & 0 & (B_A{\otimes} I_{N_{\hspace{-0.5mm}A}})C_{BA}\\
B_{AB}(C_A{\otimes} I_{N_{\hspace{-0.5mm}A}}) & A_{AB} & 0 & 0\\
0 & (B_B{\otimes} I_{N_{\hspace{-0.5mm}B}}) C_{AB} & A_B{\otimes} I_{N_{\hspace{-0.5mm}B}} & 0\\
0 & 0 & {B}_{BA}(C_B{\otimes} I_{N_{\hspace{-0.5mm}B}}) & A_{BA}
\end{array}\hspace{-1mm}\right]\hspace{-1mm},
%\left[\begin{array}{c}{X}_A\\{X}_{AB}\\{X}_B\\{X}_{BA}\end{array}\right]
\end{equation}
where matrices $A_A\in\R^{n{\times}n}$, $B_A\in\R^{n{\times}1}$, $C_A\in\R^{1{\times}n}$ are associated with the linearization of $H_A$; and matrices $A_{AB}\in\R^{(N_A+2N_B){\times}(N_A+2N_B)}$, $B_{AB}\in\R^{(N_A+2N_B){\times}N_A}$, $C_{AB}\in\R^{N_B{\times}(N_A+2N_B)}$ are the linearization matrices of the transceiver $\text{tx/rx}_{A\rightarrow B}$. For the transceiver, the linearization matrices are of the form:
\begin{equation*}
A_{AB}=\left[\begin{array}{cc}L _{AB}& \begin{array}{c}0\\0\end{array}\\ \begin{array}{cc}0&0\end{array}&0 \end{array}\right]+\left[\begin{array}{ccc}\partial\Gamma_x & 0 & 0\\ 0 &\partial\Phi_x & \partial\Phi_R\\ 0 & \partial\Psi_x & \partial\Psi_R \end{array}\right],\vspace{-2mm}
\end{equation*}
and with
\begin{equation*} 
B_{AB}= \left[\begin{array}{ccc} \partial \Gamma_u & 0_{N_B\times N_A} & 0_{N_B\times N_A} \end{array} \right]^T,
\end{equation*}
\begin{equation*}
C_{AB}=\left[\begin{array}{ccc}0_{N_B\times N_A} & 0_{N_B\times N_B} & I_{N_B}\end{array}\right],
\end{equation*}
where due to the structure of the steady state, $\partial\Gamma_x=\partial\gamma_X I_{N_A}$ with $\partial\gamma_X\triangleq\frac{\partial \gamma_X^i}{\partial X_A^i}|_{\tilde{x}_A}$, and similarly the matrices $\partial\Phi_x$, $\partial\Phi_R$, $\partial\Psi_x$, $\partial\Psi_R$, and $\partial\Gamma_u$, are diagonal with constants $\partial\phi_x$, $\partial\phi_R$, $\partial\psi_x$, $\partial\psi_R$ and $\partial\gamma_u$, respectively. The matrix $L_{AB}$ is the Laplacian matrix of the network when labeling first the nodes of type $A$.\\
 
Due to the monotonicity property of the network proved in Lemma \ref{lemma:monotonicityNetwork}, the proof follows as discussed above, and in a similar way to \cite[Proof of Theorem 2]{arcak13}. We write \eqref{eq:bigjacobianFULL} as a unitary positive feedback system: $\mathcal{A}+\mathcal{B}\mathcal{C}$ where $\mathcal{C}{=}\left[0\  0\ 0\ C_{AB}\right]$, $B{=}\left[B_A{\otimes}I_{N_{\hspace{-0.5mm}A}}^T\ 0\ 0\ 0\right]^T$, and $\mathcal{A}$ is the block triangular matrix defined in \eqref{eq:bigjacobianFULL} except for the block $(B_A{\otimes} I_{N_{\hspace{-0.5mm}A}})C_{\hspace{-0.5mm}BA}$, which is replaced by $0\in\R^{nN_A{\times}(N_B+2N_A)}$. Then, since the network is monotone with respect to the same input and output cones, we conclude stability from $-(I+\mathcal{C}\mathcal{A}^{-1}\mathcal{B})$. First note that:\par\nobreak
\vspace{-2.8mm}\hspace{-1mm}{\scriptsize
\begin{eqnarray*}
\mathcal{C}\mathcal{A}^{-1}\mathcal{B}=&\\
&\hspace{-1.45cm}={-}C_{BA}A_{BA}^{-1}B_{BA}(C_B A_B^{-1} B_B\hspace{-0.5mm}\otimes\hspace{-0.5mm} I_{N_{\hspace{-0.5mm}B}})C_{AB}A_{AB}^{-1}B_{AB}(C_A A_A^{-1} B_A\hspace{-0.5mm}\otimes\hspace{-0.5mm} I_{N_{\hspace{-0.5mm}A}})\\
&\hspace{-2.15cm}={-}T'_{A}(T_{BA}(\tilde{z}_B))T'_{B}(T_{AB}(\tilde{z}_A))(C_{BA}A_{BA}^{-1}B_{BA})(C_{AB}A_{AB}^{-1}B_{AB}),%\hspace{-0.2mm}
\end{eqnarray*}}\par\nobreak
\noindent where the second equality follows from a derivation similar to \cite{arcak13} where $T_k'(\tilde{z})=-C^{\tilde{z}}_k(A^{\tilde{z}}_k)^{-1}B^{\tilde{z}}_k$ is the static input-output map for each block at steady-state $\tilde{z}$, and $C_j$, $A_j$, $B_j$ are the linearization matrices of each block at $\tilde{z}$, we drop the superscripts $\tilde{z}$ to simplify the notation. Assumptions \ref{ass:differentiabilityTR} and \ref{ass:differentiability} guarantee that $A_j^{-1}$ exists and that $\mathcal{A}$ is nonsingular.\\
For the final step, we use the equitability assumption on the partition defined by the classes $O_A$ and $O_B$ to derive the largest eigenvalue of the matrix $(C_{BA}A_{BA}^{-1}B_{BA}C_{AB}A_{AB}^{-1}B_{AB})\in\R^{N_B\times N_B}$, and therefore the stability of the matrix $-(I+\mathcal{C}\mathcal{A}^{-1}\mathcal{B})$.
\vspace{2mm}\begin{claim}
The largest eigenvalue of the matrix $(C_{BA}A_{BA}^{-1}B_{BA}C_{AB}A_{AB}^{-1}B_{AB})$ is given by $(\overline{C}_{BA}\overline{A}_{BA}^{-1}\overline{B}_{BA}\overline{C}_{AB}\overline{A}_{AB}^{-1}\overline{B}_{AB})$ with eigenvector $\mathbf{1}_{N_A}$, where
\begin{equation*}
\overline{A}_{AB}=\left[\begin{array}{cc}\overline{L}_{AB} & \begin{array}{c}0\\0\end{array}\\ \begin{array}{cc}0&0\end{array}&0 \end{array}\right]+\left[\begin{array}{ccc}\partial\gamma_x & 0 & 0\\ 0 &\partial\phi_x & \partial\phi_R\\ 0 & \partial\psi_x & \partial\psi_R \end{array}\right],\vspace{-2mm}
\end{equation*}
and with
\begin{equation*}
\overline{B}_{AB}=\left[\begin{array}{ccc} \partial\gamma_u & 0 & 0 \end{array} \right]^T,\ \ \ 
\overline{C}_{AB}=\left[\begin{array}{ccc}0 & 0 & 1\end{array}\right],
\end{equation*}
where $\overline{A}_{AB}\in\R^{3\times 3}$, $\overline{L}_{AB}\in\R^{2\times 2}$ is the quotient Laplacian, $\overline{C}_{AB}\in\R^{1\times 3}$, and $\overline{B}_{AB}\in\R^{3\times 1}$; and by appropriate change of subscripts the same follows for the matrices $\overline{A}_{AB}$, $\overline{B}_{AB}$ and $\overline{C}_{AB}$.\hfill$\blacksquare$
\end{claim}\vspace{2mm}
The theorem follows from this claim because $T'_{AB}(\tilde{z}_A)={-}\overline{C}_{AB}\overline{A}_{AB}^{-1}\overline{B}_{AB}$, and thus the largest eigenvalue of $\mathcal{C}\mathcal{A}^{-1}\mathcal{B}$ is given by $T'_{A}(T_{BA}(\tilde{z}_B))T'_{B}(T_{AB}(\tilde{z}_A))T'_{AB}(\tilde{z}_A)T'_{BA}(\tilde{z}_B)$. Therefore, when inequality \eqref{eq:stabilityCondFULL} holds the matrix $-(I+\mathcal{C}\mathcal{A}^{-1}\mathcal{B})$ is Hurwitz and the steady-state is asymptotically stable. If the condition \eqref{eq:existenceCond} holds, $-(I+\mathcal{C}\mathcal{A}^{-1}\mathcal{B})$ has a positive eigenvalue and the steady-state is unstable.\\
\textit{Proof of Claim}: 
First note that due to equitability of the compartmental network, we can construct matrices $Q_{AB}\in\R^{(N_A+2N_B)\times 3}$ where
\begin{eqnarray*}
Q_{AB}=\hspace{-3mm}&\left[\begin{array}{ccc} 1\ ...\ 1 & 0\ ...\ 0 & 0\ ...\ 0\\ 0\ ...\ 0 & 1\ ...\ 1 & 0\ ...\ 0\\ 0\ ...\ 0 & 0\ ...\ 0 & 1\ ...\ 1\vspace{-1mm}
\end{array}\right]^T,&\\
&\hspace{-2mm}\underbrace{\hphantom{1\ ...\ 1}}_{\times N_A}\hspace{3mm} \underbrace{\hphantom{1\ ...\ 1}}_{\times N_B} \hspace{3mm}\underbrace{\hphantom{1\ ...\ 1}}_{\times N_B}&
\end{eqnarray*}
and similarly $Q_{BA}\,{\in}\,\R^{(N_B+2N_A)\times 3}$ with appropriate dimensions.
Therefore, due to equitability $L_{AB}Q_{AB}{=}Q_{AB}\overline{L}_{AB}$ and $L_{BA}Q_{BA}{=}Q_{BA}\overline{L}_{BA}$. Let $P{:=}[Q\  R]$ where $R$ is a matrix in $\R^{(N_A+2N_B)\times(N_A+2N_B-3)}$ (or $R\in\R^{(N_B+2N_A)\times(N_B+2N_A-3)}$) such that its columns, together with those of $Q$, from a basis for $\R^{N_A+2N_B}$ (or $\R^{N_B+2N_A}$). We conclude that, there exist matrices $N$ and $M$ such that
\begin{equation}
P_{AB}^{-1}A_{AB}P_{AB}=\left[\begin{array}{cc}\overline{A}_{AB}& N\\ 0 & M\end{array}\right],
\end{equation}
and similarly for $A_{BA}$. Therefore,\par\nobreak\vspace{-3mm}
{\small \begin{align*}
C_{AB}A^{-1}_{AB}&B_{AB}\mathbf{1}_{N_A} =\\
&=(C_{AB}P_{AB})(P_{AB}^{-1}A_{AB}P_{AB})^{-1}(P_{AB}^{-1}B_{AB}\mathbf{1}_{N_A})\\
&= \left[\begin{array}{cc}\overline{C}_{AB}\mathbf{1}_{N_B} & S\end{array}\right] \left[\begin{array}{cc} \overline{A}_{AB}^{-1} &  U\\ 0 & V \end{array}\right] \left[\begin{array}{c}\overline{B}_{AB}\\0\end{array}\right],\\
&= \overline{C}_{AB}\overline{A}_{AB}^{-1}\overline{B}_{AB}\mathbf{1}_{N_B}.
\end{align*}}\par\nobreak
\noindent for some matrices $S$, $U$, and $V$ with appropriate dimensions. This implies that\par\nobreak\vspace{-3mm}
{\small \begin{align*}
C_{BA}A^{-1}_{BA}B_{BA}&C_{AB}A^{-1}_{AB}B_{AB}\mathbf{1}_{N_A} =\\
&=  (\overline{C}_{AB}\overline{A}_{AB}^{-1}\overline{B}_{AB})C_{BA}A^{-1}_{BA}B_{BA}\mathbf{1}_{N_B}\\
&= (\overline{C}_{BA}\overline{A}_{BA}^{-1}\overline{B}_{BA}\overline{C}_{AB}\overline{A}_{AB}^{-1}\overline{B}_{AB})\mathbf{1}_{N_A}
\end{align*}}\par\nobreak
\noindent\textit{i.e.}, $\overline{C}_{BA}\overline{A}_{BA}^{-1}\overline{B}_{BA}\overline{C}_{AB}\overline{A}_{AB}^{-1}\overline{B}_{AB}=T'_{AB}(\tilde{z}_A)T'_{BA}(\tilde{z}_B)$ is an eigenvalue of $C_{BA}A^{-1}_{BA}B_{BA}C_{AB}A^{-1}_{AB}B_{AB}$ with associated eigenvector $\mathbf{1}_{N_A}$. Note that this eigenvalue is positive since the static input/output maps of the transceivers have positive slope. Finally, we need to show that this is the largest eigenvalue. Note that due to Assumption \ref{ass:monotonicityTR}, the transceivers' input/output maps $T_{AB}^{\text{tx/rx}}(\mathbf{\tilde{z}}_A){=}{-}C_{AB}A^{-1}_{AB}B_{AB}$ and $T_{BA}^{\text{tx/rx}}(\mathbf{\tilde{z}}_B){=}{-}C_{BA}A^{-1}_{BA}B_{BA}$ are nonnegative matrices \cite{angeli03}, and thus so is $T_{AB}^{\text{tx/rx}}(\mathbf{\tilde{z}}_A)T_{BA}^{\text{tx/rx}}(\mathbf{\tilde{z}}_B)$, with no zero rows. This concludes the proof of the claim since, by the Perron-Frobenius Theorem \cite{berman94}, %[Proof of 2.10]
 the eigenvalue with associated positive eigenvector $\mathbf{1}_{N_A}$, must be the largest positive eigenvalue.\hfill$\blacksquare$

%\section*{Appendix}
\begin{table}[ht]
\setlength{\tabcolsep}{2pt}
\centering
\caption{Parameters used in simulations}
\label{table:parameters}
\scriptsize
\begin{tabular}{| c | c | c | c |} %V_{\text{P}_{\textit{LuxI}}}   %V_{\text{P}_{\textit{LtetO-}1}}
\hline
Parameter & Description & Value & Units\\\hline
$k_{on}$& binding rate between {LuxR} and {AHL} &
$1e{9}$ & s$^{\mbox{-}1}$M$^{\mbox{-}1}$\\\hline 
$k_{off}$ & dissociation rate between {LuxR} and {AHL} & $50$
& s$^{\mbox{-}1}$\\\hline 
$p_{R_i}$ & constitutive concentration of total {LuxR} & variable & M\\\hline
$d_{12}$ & diffusion rate of {AHL} & variable & s$^{\mbox{-}1}$\\\hline
$V_{\text{P}_{\text{LuxI}}}$ & velocity rate of promoter  P$_{\text{LuxI}}$ & $0.26$ & s$^{\mbox{-}1}$\\\hline
$N_{\text{P}_{\text{LuxI}}}$ & copy number of promoter P$_{\text{LuxI}}$ & $5$ & $1$\\\hline
$C$ & concentration of a single molecule in a cell & $1.5e{\mbox{-}9}$ & M\\\hline 
$K_{RA}$ & dissociation constant between p$_{R_A}$ and P$_{\text{LuxI}}$ & $1.5e{\mbox{-}9}$ & M\\\hline 
$n_{RA}$ & cooperativity & $2$ & $1$\\\hline
$\ell_{\text{P}_{\text{LuxI}}}$ & leakage of promoter P$_{\text{LuxI}}$ & $1/167$ & $1$\\\hline
$V_{\text{P}_{\text{LtetO-}1}}$ & velocity rate of promoter P$_{\text{LtetO-}1}$ & $0.3$ & s$^{\mbox{-}1}$\\\hline
$N_{\text{P}_{\text{LtetO-}1}}$ & copy number of promoter P$_{\text{LtetO-}1}$ & $5$ & $1$\\\hline
$K_T$ & dissociation constant between {TetR} and P$_{\text{LtetO-}1}$ & $1.786e{\mbox{-}10}$ & M\\\hline 
$n_T$ & cooperativity & $2$ & $1$\\\hline
$\ell_{\text{P}_{\text{LtetO-}1}}$ & leakage of promoter P$_{\text{P}_{\text{LtetO-}1}}$ & $1/5050$ & $1$\\\hline
$\gamma_{A}$ & rate of degradation of {AHL} & $7.70e{\mbox{-}4}$ &
s$^{\mbox{-}1}$\\\hline 
$\gamma_{m_T}$ & degradation constant of mRNA {TetR} & $5.78e{\mbox{-}3}$ &
s$^{\mbox{-}1}$\\\hline 
$\gamma_T$ & degradation constant of {TetR} & $2.89e{\mbox{-}4}$ &
s$^{\mbox{-}1}$\\\hline 
$\gamma_{m_I}$ & degradation constant of mRNA {LuxI}/BviI & $5.78e{\mbox{-}3}$ &
s$^{\mbox{-}1}$\\\hline 
$\gamma_{I}$ & degradation constant of {LuxI}/BviI & $1.16e{\mbox{-}3}$ &
s$^{\mbox{-}1}$\\\hline 
$\epsilon_{T}$ & translation rate \textit{tetR} & $6.224e{\mbox{-}6}$ &
s$^{\mbox{-}1}$\\\hline 
$\epsilon_{I}$ & translation rate \textit{luxI}/\textit{bviI} & $2.655e{\mbox{-}5}$ &
s$^{\mbox{-}1}$\\\hline 
$\nu$ & generation rate of {AHL} & $0.0135$ & s$^{\mbox{-}1}$\\\hline
\end{tabular}
\end{table}

\bibliographystyle{ieeetr}
\bibliography{doc}

\end{document}